\documentclass[a4paper,oneside]{amsart}
\usepackage[T1]{fontenc}
\usepackage[utf8]{inputenc}
\usepackage[active]{srcltx}
\usepackage{units}
\usepackage{amsthm}
\usepackage{amstext}
\usepackage{amssymb}
\usepackage{esint}
\usepackage{xargs}[2008/03/08]

\makeatletter

\numberwithin{equation}{section}
\numberwithin{figure}{section}
\theoremstyle{plain}
\newtheorem{thm}{\protect\theoremname}
  \theoremstyle{plain}
  \newtheorem{cor}[thm]{\protect\corollaryname}
  \theoremstyle{remark}
  \newtheorem*{rem*}{\protect\remarkname}

\makeatother

  \providecommand{\corollaryname}{Corollary}
  \providecommand{\remarkname}{Remark}
\providecommand{\theoremname}{Theorem}

\begin{document}

\title{Daehee numbers and polynomials}

\author{Dae San Kim and Taekyun Kim }
\begin{abstract}
We consider the Witt-type formula for Daehee numbrers and polynomials
and investigate some properties of those numbers and polynomials.
In particular, Daehee numbers are closely related to higher-order
Bernoulli numbers and Bernoulli numbers of the second kind.

\newcommandx\ud[1][usedefault, addprefix=\global, 1=]{\textnormal{UD}\left(#1\right)}

\end{abstract}
\maketitle

\section{Introduction}

As is known, the $n$-th Daehee polynomials are defined by the generating
function to be
\begin{equation}
\left(\frac{\log\left(1+t\right)}{t}\right)\left(1+t\right)^{x}=\sum_{n=0}^{\infty}D_{n}\left(x\right)\frac{t^{n}}{n!},\:\left(\textrm{see\:[5,6,8,9,10,11]}\right).\label{eq:1}
\end{equation}

In the special case, $x=0$, $D_{n}=D_{n}\left(0\right)$ are called
the Daehee numbers.

Throughout this paper, $\mathbb{Z}_{p}$, $\mathbb{Q}_{p}$ and $\mathbb{C}_{p}$
will denote the rings of $p$-adic integers, the fields of $p$-adic
numbers and the completion of algebraic closure of $\mathbb{Q}_{p}$.
The $p-$adic norm $\left|\cdot\right|_{p}$ is normalized by
$\left|p\right|_{p}=\nicefrac{1}{p}$. Let $\ud[\mathbb{Z}_{p}]$
be the space of uniformly differentiable functions on $\mathbb{Z}_{p}$.
For $f\in\ud[\mathbb{Z}_{p}]$, the $p$-adic invariant integral on
$\mathbb{Z}_{p}$ is defined by
\begin{equation}
I\left(f\right)=\int_{\mathbb{Z}_{p}}f\left(x\right)d\mu_{0}\left(x\right)=\lim_{n\rightarrow\infty}\frac{1}{p^{n}}\sum_{x=0}^{p^{n}-1}f\left(x\right),\:\left(\textrm{see }[6]\right).\label{eq:2}
\end{equation}

Let $f_{1}$ be the translation of $f$ with $f_{1}\left(x\right)=f\left(x+1\right).$
Then, by (\ref{eq:2}), we get
\begin{equation}
I\left(f_{1}\right)=I\left(f\right)+f^{\prime}\left(0\right),\:\textrm{where }f^{\prime}\left(0\right)=\left.\frac{df\left(x\right)}{dx}\right|_{x=0}.\label{eq:3}
\end{equation}

As is known, the Stirling number of the first kind is defined by
\begin{equation}
\left(x\right)_{n}=x\left(x-1\right)\cdots\left(x-n+1\right)=\sum_{l=0}^{n}S_{1}\left(n,l\right)x^{l},\label{eq:4}
\end{equation}
 and the Stirling number of the second kind is given by the generating
function to be
\begin{equation}
\left(e^{t}-1\right)^{m}=m!\sum_{l=m}^{\infty}S_{2}\left(l,m\right)\frac{t^{l}}{l!},\:\left(\textrm{see [2,3,4]}\right).\label{eq:5}
\end{equation}

For $\alpha\in\mathbb{Z}$, the Bernoulli polynomials of order
$\alpha$ are defined by the generating function to be
\begin{equation}
\left(\frac{t}{e^{t}-1}\right)^{\alpha}e^{xt}=\sum_{n=0}^{\infty}B_{n}^{\left(\alpha\right)}\left(x\right)\frac{t^{n}}{n!},\:\left(\textrm{see }[1,2,8]\right).\label{eq:6}
\end{equation}

When $x=0$, $B_{n}^{\left(\alpha\right)}=B_{n}^{\left(\alpha\right)}\left(0\right)$
are called the Bernoulli numbers of order $\alpha$.

In this paper, we give a $p$-adic integral representation of Daehee
numbers and polynomials, which are called the Witt-type formula for Daehee
numbers and polynomials. From our integral representation, we can
derive some interesting properties related to Daehee numbers and polynomials.

\section {Witt-type formula for Daehee numbers and polynomials}

First, we consider the following integral representation associated
with falling factorial sequences :
\begin{equation}
\int_{\mathbb{Z}_{p}}\left(x\right)_{n}d\mu_{0}\left(x\right),\:\textrm{ where }n\in\mathbb{Z}_{+}=\mathbb{N}\cup\left\{ 0\right\} .\label{eq:7}
\end{equation}

By (\ref{eq:7}), we get
\begin{eqnarray}
\sum_{n=0}^{\infty}\int_{\mathbb{Z}_{p}}\left(x\right)_{n}d\mu_{0}\left(x\right)\frac{t^{n}}{n!} & = & \int_{\mathbb{Z}_{p}}\sum_{n=0}^{\infty}\dbinom{x}{n}t^{n}d\mu_{0}\left(x\right)\label{eq:8}\\
 & = & \int_{\mathbb{Z}_{p}}\left(1+t\right)^{x}d\mu_{0}\left(x\right),\nonumber
\end{eqnarray}
where $t\in\mathbb{C}_{p}$ with $\left|t\right|_{p}<p^{-\frac{1}{p-1}}.$

For $t\in\mathbb{C}_{p}$ with $\left|t\right|_{p}<p^{-\frac{1}{p-1}}$,
let us take $f\left(x\right)=\left(1+t\right)^{x}.$ Then, from (\ref{eq:3}),
we have
\begin{equation}
\int_{\mathbb{Z}_{p}}\left(1+t\right)^{x}d\mu_{0}\left(x\right)=\frac{\log\left(1+t\right)}{t}.\label{eq:9}
\end{equation}

By (\ref{eq:1}) and (\ref{eq:9}), we see that
\begin{eqnarray}
\sum_{n=0}^{\infty}D_{n}\frac{t^{n}}{n!} & = & \frac{\log\left(1+t\right)}{t}\label{eq:10}\\
 & = & \int_{\mathbb{Z}_{p}}\left(1+t\right)^{x}d\mu_{0}\left(x\right)\nonumber \\
 & = & \sum_{n=0}^{\infty}\int_{\mathbb{Z}_{p}}\left(x\right)_{n}d\mu_{0}\left(x\right)\frac{t^{n}}{n!}.\nonumber
\end{eqnarray}

Therefore, by (\ref{eq:10}), we obtain the following theorem.
\begin{thm}
\label{thm:1}For $n\ge0$, we have
\[
\int_{\mathbb{Z}_{p}}\left(x\right)_{n}d\mu_{0}\left(x\right)=D_{n}.
\]

\end{thm}

For $n\in\mathbb{Z}$, it is known that
\begin{equation}
\left(\frac{t}{\log\left(1+t\right)}\right)^{n}\left(1+t\right)^{x-1}=\sum_{k=0}^{\infty}B_{k}^{\left(k-n+1\right)}\left(x\right)\frac{t^k}{k!},\:\left(\textrm{see [2,3,4]}\right).\label{eq:11}
\end{equation}

Thus, by (\ref{eq:11}), we get
\begin{equation}
D_{k}=\int_{\mathbb{Z}_{p}}\left(x\right)_{k}d\mu_{0}\left(x\right)=B_{k}^{\left(k+2\right)}\left(1\right),\quad\left(k\ge0\right),\label{eq:12}
\end{equation}
where $B_{k}^{\left(n\right)}\left(x\right)$ are the Bernoulli
polynomials of order $n$.

In the special case, $x=0$, $B_{k}^{\left(n\right)}=B_{k}^{\left(n\right)}\left(0\right)$
are called the $n$-th Bernoulli numbers of order $n$.

From (\ref{eq:10}), we note that
\begin{eqnarray}
\left(1+t\right)^{x}\int_{\mathbb{Z}_{p}}\left(1+t\right)^{y}d\mu_{0}\left(y\right) & = & \left(\frac{\log\left(1+t\right)}{t}\right)\left(1+t\right)^{x}\label{eq:13}\\
 & = & \sum_{n=0}^{\infty}D_{n}\left(x\right)\frac{t^{n}}{n!}.\nonumber
\end{eqnarray}
 Thus, by (\ref{eq:13}), we get
\begin{equation}
\int_{\mathbb{Z}_{p}}\left(x+y\right)_{n}d\mu_{0}\left(y\right)=D_{n}\left(x\right),\quad\left(n\ge0\right),\label{eq:14}
\end{equation}
 and, from (\ref{eq:11}), we have
\begin{equation}
D_{n}\left(x\right)=B_{n}^{\left(n+2\right)}\left(x+1\right).\label{eq:15}
\end{equation}

Therefore, by (\ref{eq:14}) and (\ref{eq:15}), we obtain the following
theorem.
\begin{thm}
\label{thm:2}For $n\ge0$, we have
\[
D_{n}\left(x\right)=\int_{\mathbb{Z}_{p}}\left(x+y\right)_{n}d\mu_{0}\left(y\right),
\]
 and
\[
D_{n}\left(x\right)=B_{n}^{\left(n+2\right)}\left(x+1\right).
\]

\end{thm}

By Theorem \ref{thm:1}, we easily see that
\begin{equation}
D_{n}=\sum_{l=0}^{n}S_{1}\left(n,l\right)B_{l},\label{eq:16}
\end{equation}
where $B_{l}$ are the ordinary Bernoulli numbers.

From Theorem \ref{thm:2}, we have
\begin{eqnarray}
D_{n}\left(x\right) & = & \int_{\mathbb{Z}_{p}}\left(x+y\right)_{n}d\mu_{0}\left(y\right)\label{eq:17}\\
 & = & \sum_{l=0}^{n}S_{1}\left(n,l\right)B_{l}\left(x\right),\nonumber
\end{eqnarray}
where $B_{l}\left(x\right)$ are the Bernoulli polynomials
defined by generating function to be
\[
\frac{t}{e^{t}-1}e^{xt}=\sum_{n=0}^{\infty}B_{n}\left(x\right)\frac{t^{n}}{n!}.
\]

Therefore, by (\ref{eq:16}) and (\ref{eq:17}), we obtain the following
corollary.
\begin{cor}
\label{cor:3}For $n\ge0$, we have
\[
D_{n}\left(x\right)=\sum_{l=0}^{n}S_{1}\left(n,l\right)B_{l}\left(x\right).
\]

\end{cor}

In (\ref{eq:10}), we have
\begin{eqnarray}
\frac{t}{e^{t}-1} & = & \sum_{n=0}^{\infty}D_{n}\frac{1}{n!}\left(e^{t}-1\right)^{n}\label{eq:18}\\
 & = & \sum_{n=0}^{\infty}D_{n}\frac{1}{n!}n!\sum_{m=n}^{\infty}S_{2}\left(m,n\right)\frac{t^{m}}{m!}\nonumber \\
 & = & \sum_{m=0}^{\infty}\left(\sum_{n=0}^{m}D_{n}S_{2}\left(m,n\right)\right)\frac{t^{m}}{m!}\nonumber
\end{eqnarray}
and
\begin{equation}
\frac{t}{e^{t}-1}=\sum_{m=0}^{\infty}B_{m}\frac{t^{m}}{m!}.\label{eq:19}
\end{equation}

Therefore, by (\ref{eq:18}) and (\ref{eq:19}), we obtain the following
theorem.
\begin{thm}
\label{thm:4}For $m\ge0$, we have
\[
B_{m}=\sum_{n=0}^{m}D_{n}S_{2}\left(m,n\right).
\]

In particular,
\[
\int_{\mathbb{Z}_{p}}x^{m}d\mu_{0}\left(x\right)=\sum_{n=0}^{m}D_{n}S_{2}\left(m,n\right).
\]
\end{thm}
\begin{rem*}
For $m\ge0$, by (\ref{eq:17}), we have
\[
\int_{\mathbb{Z}_{p}}\left(x+y\right)^{m}d\mu_{0}\left(y\right)=\sum_{n=0}^{m}D_{n}\left(x\right)S_{2}\left(m,n\right).
\]

\end{rem*}

For $n\in\mathbb{Z}_{\ge0}$, the rising factorial sequence is defined
by
\begin{equation}
x^{\left(n\right)}=x\left(x+1\right)\cdots\left(x+n-1\right).\label{eq:20}
\end{equation}

Let us define the Daehee numbers of the second kind as follows :
\begin{equation}
\widehat{D}_{n}=\int_{\mathbb{Z}_{p}}\left(-x\right)_{n}d\mu_{0}\left(x\right),\:\left(n\in\mathbb{Z}_{\ge0}\right).\label{eq:21}
\end{equation}

By (\ref{eq:21}), we get
\begin{equation}
x^{\left(n\right)}=\left(-1\right)^{n}\left(-x\right)_{n}=\sum_{l=0}^{n}S_{1}\left(n,l\right)\left(-1\right)^{n-l}x^{l}.\label{eq:22}
\end{equation}

From (\ref{eq:21}) and (\ref{eq:22}), we have
\begin{eqnarray}
\widehat{D}_{n} & = & \int_{\mathbb{Z}_{p}}\left(-x\right)_{n}d\mu_{0}\left(x\right)=\int_{\mathbb{Z}_{p}}x^{\left(n\right)}\left(-1\right)^{n}d\mu_{0}\left(x\right)\label{eq:23}\\
 & = & \sum_{l=0}^{n}S_{1}\left(n,l\right)\left(-1\right)^{l}B_{l}.\nonumber
\end{eqnarray}

Therefore, by (\ref{eq:23}), we obtain the following theorem.
\begin{thm}
\label{thm:5}For $n\ge0$, we have
\[
\widehat{D}_{n}=\sum_{l=0}^{n}S_{1}\left(n,l\right)\left(-1\right)^{l}B_{l}.
\]

\end{thm}

Let us consider the generating function of the Daehee numbers of the
second kind as follows :

\begin{eqnarray}
\sum_{n=0}^{\infty}\widehat{D}_{n}\frac{t^{n}}{n!} & = & \sum_{n=0}^{\infty}\int_{\mathbb{Z}_{p}}\left(-x\right)_{n}d\mu_{0}\left(x\right)\frac{t^{n}}{n!}\label{eq:24}\\
 & = & \int_{\mathbb{Z}_{p}}\sum_{n=0}^{\infty}\dbinom{-x}{n}t^{n}d\mu_{0}\left(x\right)\nonumber \\
 & = & \int_{\mathbb{Z}_{p}}\left(1+t\right)^{-x}d\mu_{0}\left(x\right).\nonumber
\end{eqnarray}

From (\ref{eq:3}), we can derive the following equation :

\begin{equation}
\int_{\mathbb{Z}_{p}}\left(1+t\right)^{-x}d\mu_{0}\left(x\right)=\frac{\left(1+t\right)\log\left(1+t\right)}{t},\label{eq:25}
\end{equation}
where $\left|t\right|_{p}<p^{-\frac{1}{p}}.$

By (\ref{eq:24}) and (\ref{eq:25}), we get
\begin{eqnarray}
\frac{1}{t}\left(1+t\right)\log\left(1+t\right) & = & \int_{\mathbb{Z}_{p}}\left(1+t\right)^{-x}d\mu_{0}\left(x\right)\label{eq:26}\\
 & = & \sum_{n=0}^{\infty}\widehat{D}_{n}\frac{t^{n}}{n!}.\nonumber
\end{eqnarray}

Let us consider the Daehee polynomials of the second kind as follows
:

\begin{equation}
\frac{\left(1+t\right)\log\left(1+t\right)}{t}\frac{1}{\left(1+t\right)^{x}}=\sum_{n=0}^{\infty}\widehat{D}_{n}\left(x\right)\frac{t^{n}}{n!}.\label{eq:27}
\end{equation}

Then, by (\ref{eq:27}), we get
\begin{equation}
\int_{\mathbb{Z}_{p}}\left(1+t\right)^{-x-y}d\mu_{0}\left(y\right)=\sum_{n=0}^{\infty}\hat{D}_{n}\left(x\right)\frac{t^{n}}{n!}.\label{eq:28}
\end{equation}

From (\ref{eq:28}), we get
\begin{eqnarray}
\widehat{D}_{n}\left(x\right) & = & \int_{\mathbb{Z}_{p}}\left(-x-y\right)_{n}d\mu_{0}\left(y\right),\quad\left(n\ge0\right)\label{eq:29}\\
 & = & \sum_{l=0}^{n}\left(-1\right)^{l}S_{1}\left(n,l\right)B_{l}\left(x\right).\nonumber
\end{eqnarray}

Therefore, by (\ref{eq:29}), we obtain the following theorem.
\begin{thm}
\label{thm:6}For $n\ge0$, we have
\[
\widehat{D}_{n}\left(x\right)=\int_{\mathbb{Z}_{p}}\left(-x-y\right)_{n}d\mu_{0}\left(y\right)=\sum_{l=0}^{n}\left(-1\right)^{l}S_{1}\left(n,l\right)B_{l}\left(x\right).
\]

\end{thm}

From (\ref{eq:27}) and (\ref{eq:28}), we have
\begin{eqnarray}
\left(\frac{t}{e^{t}-1}\right)e^{\left(1-x\right)t} & = & \sum_{n=0}^{\infty}\widehat{D}_{n}\left(x\right)\frac{1}{n!}\left(e^{t}-1\right)^{n}\label{eq:30}\\
 & = & \sum_{n=0}^{\infty}\widehat{D}_{n}\left(x\right)\frac{1}{n!}n!\sum_{m=n}^{\infty}S_{2}\left(m,n\right)\frac{t^{m}}{m!}\nonumber \\
 & = & \sum_{m=0}^{\infty}\left(\sum_{n=0}^{m}\widehat{D}_{n}\left(x\right)S_{2}\left(m,n\right)\right)\frac{t^{n}}{m!},\nonumber
\end{eqnarray}
 and
\begin{eqnarray}
\int_{\mathbb{Z}_{p}}e^{-\left(x+y\right)t}d\mu_{0}\left(y\right) & = & \sum_{n=0}^{\infty}\widehat{D}_{n}\left(x\right)\frac{\left(e^{t}-1\right)^{n}}{n!}\label{eq:31}\\
 & = & \sum_{m=0}^{\infty}\left(\sum_{n=0}^{m}\widehat{D}_{n}\left(x\right)S_{2}\left(m,n\right)\right)\frac{t^{m}}{m!}.\nonumber
\end{eqnarray}

Therefore, by (\ref{eq:30}) and (\ref{eq:31}), we obtain the follwoing
theorem.
\begin{thm}
\label{thm:7}For $m\ge0$, we have
\begin{eqnarray*}
B_{m}\left(1-x\right) & = & \left(-1\right)^{m}\int_{\mathbb{Z}_{p}}\left(x+y\right)^{m}d\mu_{0}\left(y\right)\\
 & = & \sum_{n=0}^{m}\widehat{D}_{n}\left(x\right)S_{2}\left(m,n\right).
\end{eqnarray*}

In particular,
\[
B_{m}\left(1-x\right)=\left(-1\right)^{m}B_{m}\left(x\right)=\sum_{n=0}^{m}\widehat{D}_{m}\left(x\right)S_{2}\left(m,n\right).
\]
\end{thm}
\begin{rem*}
By (\ref{eq:11}), (\ref{eq:26}) and (\ref{eq:27}), we see that
\[
\widehat{D}_{n}=B_{n}^{\left(n+2\right)}\left(2\right),\quad\widehat{D}_{n}\left(x\right)=B_{n}^{\left(n+2\right)}\left(2-x\right).
\]

\end{rem*}

From Theorem \ref{thm:1} and (\ref{eq:21}), we have
\begin{eqnarray}
\left(-1\right)^{n}\frac{D_{n}}{n!} & = & \left(-1\right)^{n}\int_{\mathbb{Z}_{p}}\dbinom{x}{n}d\mu_{0}\left(x\right)\label{eq:32}\\
 & = & \int_{\mathbb{Z}_{p}}\dbinom{-x+n-1}{n}d\mu_{0}\left(x\right)\nonumber \\
 & = & \sum_{m=0}^{n}\dbinom{n-1}{n-m}\int_{\mathbb{Z}_{p}}\dbinom{-x}{m}d\mu_{0}\left(x\right)\nonumber \\
 & = & \sum_{m=0}^{n}\dbinom{n-1}{n-m}\frac{\widehat{D}_{m}}{m!}=\sum_{m=1}^{n}\dbinom{n-1}{m-1}\frac{\widehat{D}_{m}}{m!},\nonumber
\end{eqnarray}
 and
\begin{eqnarray}
\left(-1\right)^{n}\frac{\widehat{D}_{n}}{n!} & = & \left(-1\right)^{n}\int_{\mathbb{Z}_{p}}\dbinom{-x}{n}d\mu_{0}\left(x\right)=\int_{\mathbb{Z}_{p}}\dbinom{x+n-1}{n}d\mu_{0}\left(x\right)\label{eq:33}\\
 & = & \sum_{m=0}^{n}\dbinom{n-1}{n-m}\int_{0}^{1}\dbinom{x}{m}d\mu_{0}\left(x\right)\nonumber \\
 & = & \sum_{m=0}^{n}\dbinom{n-1}{m-1}\frac{D_{m}}{m!}=\sum_{m=1}^{n}\dbinom{n-1}{m-1}\frac{D_{m}}{m!}.\nonumber
\end{eqnarray}

Therefore, by (\ref{eq:32}) and (\ref{eq:33}), we obtain the following
theorem.
\begin{thm}
\label{thm:8}For $n\in\mathbb{N},$ we have
\[
\left(-1\right)^{n}\frac{D_{n}}{n!}=\sum_{m=1}^{n}\dbinom{n-1}{m-1}\frac{\widehat{D}_{m}}{m!},
\]
and
\[
\left(-1\right)^{n}\frac{\widehat{D}_{n}}{n!}=\sum_{m=1}^{n}\dbinom{n-1}{m-1}\frac{D_{m}}{m!}.
\]
\end{thm}

$\,$

\noindent Department of Mathematics, Sogang University, Seoul 121-742, Republic of Korea

\noindent e-mail:dskim@sogang.ac.kr\\

\noindent Department of Mathematics, Kwangwoon University, Seoul 139-701, Republic of Korea

\noindent e-mail:tkkim@kw.ac.kr\\

\end{document}